\documentclass[11pt]{article}
\usepackage[colorlinks,
            linkcolor=blue,
            anchorcolor=green,
            citecolor=red
            ]{hyperref}
\usepackage{mathrsfs,amssymb,amsfonts,subfigure}
\usepackage{cite}
\usepackage{amsmath,amssymb,latexsym,color}
\usepackage{graphicx}
\usepackage[mathscr]{eucal}
\usepackage{xypic}
\newtheorem{thm}{Theorem}[section]
\newtheorem{defin}[thm]{Definition}
\newtheorem{prop}[thm]{Proposition}
\newtheorem{lemma}[thm]{Lemma}

\newtheorem{example}[thm]{Example}
\newtheorem{obser}{Observation}

\newcommand{\proof}{{\it Proof.\quad}}
\newcommand{\qed}{\hfill\Box\medskip}
\newcommand{\aut}{{\rm Aut}}
\usepackage{CJK}
\begin{document}
\begin{CJK*}{GBK}{song}
\renewcommand{\abovewithdelims}[2]{
\genfrac{[}{]}{0pt}{}{#1}{#2}}

\title{\bf The full automorphism group of the power (di)graph of a finite group}

\author{Min Feng\quad Xuanlong Ma\quad Kaishun Wang\footnote{Corresponding author. \newline {\em E-mail address:} fgmn\_1998@163.com (M. Feng), mxl881112@126.com (X. Ma), wangks@bnu.edu.cn (K. Wang).}\\
{\footnotesize   \em  Sch. Math. Sci. {\rm \&} Lab. Math. Com. Sys., Beijing Normal University, Beijing, 100875,  China} }
 \date{}
 \maketitle

\begin{abstract}
 We describe  the full automorphism group of the power (di)graph of a finite group. As an application, we solve a conjecture proposed by Doostabadi, Erfanian and   Jafarzadeh in 2013.


\medskip

\end{abstract}

\section{Introduction}
We always use $G$ to denote a finite group.
The {\em power digraph} $\overrightarrow{\mathcal P}_G$ has
 $G$ as its vertex set, where
there is an arc from $x$ to $y$ if $x\neq y$ and $y$ is a power of $x$.
The {\em  power graph}  $\mathcal P_G$ is the underlying graph of $\overrightarrow{\mathcal P}_G$, which  is obtained from
$\overrightarrow{\mathcal P}_G$ by suppressing the orientation of each arc and replacing multiple edges by one edge.
Kelarev and Quinn \cite{kel1,kel2} introduced the power digraph of a semigroup and
called it directed power graph.
Chakrabarty, Ghosh and Sen \cite{cha} defined
 power graphs of semigroups.
Recently,  power (di)graphs have been investigated by researchers, see \cite{came,mir,mog,tam,kel3}.
A detailed list of results and open questions  can be found in \cite{aba}.

In 2013,  Doostabadi, Erfanian and   Jafarzadeh asserted that the full automorphism group of the power graph of the cyclic group $Z_n$ is isomorphic to the direct product of some symmetry groups.

\medskip\noindent{\bf Conjecture} \cite{doo}  For every positive integer $n$,
\begin{equation*}
  \aut(\mathcal P_{Z_n})\cong S_{\varphi(n)+1}\times\prod_{d\in D(n)\setminus\{1,n\}} S_{\varphi(d)},
\end{equation*}
where $D(n)$ is the set of positive divisors of $n$, and $\varphi$ is the Euler's totient function.

\medskip

In fact, if $n$ is a prime power, then $\mathcal P_{Z_n}$ is a complete graph by \cite[Theorem 2.12]{cha}, which implies that $\aut(\mathcal P_{Z_n})\cong S_n$. Hence, the conjecture does not hold if $n=p^m$ for any prime $p$ and integer $m\geq 2$. The motivation of this paper is to show that this conjecture holds for the remaining case.

 In this paper we describe the full automorphism group  of  the power (di)graph  of an arbitrary finite group. As an application,   this conjecture is valid    if $n$ is not a prime power.

\section{Main results}

 Denote by $\mathcal C(G)$ the set of all cyclic subgroups of $G$. For $C\in\mathcal C(G)$, let $[C]$ denote the set of all generators of $C$. Write
$$
\mathcal C(G)=\{C_1,\ldots,C_k\}\textup{ and }[C_i]=\{[C_i]_1,\ldots,[C_i]_{s_i}\}.
$$

Define $P(G)$ as the set of permutations $\sigma$ on $\mathcal C(G)$ preserving order, inclusion and noninclusion, i.e., $|C_i^\sigma|=|C_i|$ for each $i\in\{1,\ldots,k\}$, and $C_i\subseteq C_j$ if and only if $C_i^\sigma\subseteq C_j^\sigma$. Note that $P(G)$ is a permutation group on $\mathcal C(G)$. This group induces the faithful  action on the set $G$:
\begin{equation}\label{p(G)}
G\times P(G)\longrightarrow G,\quad([C_i]_j,\sigma)\longmapsto [C_i^\sigma]_j.
\end{equation}

 For $\Omega\subseteq G$, let $S_{\Omega}$ denote the symmetric group on $\Omega$.
Since $G$ is the disjoint union of $[C_1],\ldots,[C_k]$, we get the faithful group action on the set $G$:
\begin{equation}\label{s1}
G\times\prod_{i=1}^kS_{[C_i]}\longrightarrow G,\quad([C_i]_j,(\xi_1,\ldots,\xi_k))\longmapsto ([C_i]_j)^{\xi_i}.
\end{equation}

\begin{thm}\label{mainthm1}
  Let $G$ be a finite group. Then
  \begin{eqnarray*}
  \aut(\overrightarrow{\mathcal P}_G)=(\prod_{i=1}^k S_{[C_i]})\rtimes P(G),
  \end{eqnarray*}
  where $P(G)$  and $\prod_{i=1}^k S_{[C_i]}$ act  on $G$ as in {\rm (\ref{p(G)})} and {\rm (\ref{s1})}, respectively.
\end{thm}

In the power graph $\mathcal P_G$, the {\em closed neighborhood}  of a vertex $x$, denoted by $N[x]$, is the set of its neighbors and itself. For $x,y\in G$,  define $x\equiv y$ if  $N[x]=N[y]$. Observe  that $\equiv$ is an equivalence relation. Let $\overline{x}$ denote the   equivalence class   containing $x$. Write
\begin{equation*}
\mathcal U(G)=\{\overline x\mid x\in G\}=\{\overline{u_1},\ldots,\overline{u_l}\}.
\end{equation*}
Since $G$ is the disjoint union of $\overline{u_1},\ldots,\overline{u_l}$, the following is a faithful group action on the set $G$:
\begin{equation}\label{s2}
G\times\prod_{i=1}^lS_{\overline{u_i}}\longrightarrow G,\quad(x,(\tau_1,\ldots,\tau_l))\longmapsto x^{\tau_i},\quad\textup{where }x\in\overline{u_i}.
\end{equation}

%Noting that (\ref{s1}), (\ref{p(G)}) and  (\ref{s2}) are faithful,  we may   assume that $\prod_{i=1}^kS_{[C_i]}, P(G)$ %and $\prod_{i=1}^lS_{\overline{u_i}}$ are subgroups of $S_G$.   A group $H$ is a {\em semidirect product} $H_1\rtimes %H_2$ of a normal subgroup $H_1$ by a subgroup $H_2$ if $H=H_1H_2$ and $|H_1\cap H_2|=1$.  Now we state our main result.

\begin{thm}\label{mainthm2}
  Let $G$ be a finite group. Then
  \begin{eqnarray*}
  \aut(\mathcal P_G)=(\prod_{i=1}^l S_{\overline{u_i}})\rtimes P(G),
  \end{eqnarray*}
  where $P(G)$ and $\prod_{i=1}^l S_{\overline{u_i}}$ act  on $G$ as in  {\rm (\ref{p(G)})} and {\rm (\ref{s2})}, respectively.
\end{thm}

The rest of this paper is organized as follows.   In Section 3,
   the  induced   action of $\aut ({\mathcal P}_G)$ on $\mathcal U(G)$ is discussed.
 In Section 4,  we prove  Theorems~\ref{mainthm1} and \ref{mainthm2}. In Section 5, we determine $\aut(\overrightarrow{\mathcal P}_G)$ and $\aut(\mathcal P_G)$ when $G$ is  cyclic,  elementary abelian,   dihedral or  generalized quaternion.

\section{The  induced   action of $\aut ({\mathcal P}_G)$ on $\mathcal U(G)$}

In \cite{cam},   Cameron proved that each element of $\mathcal U(G)$ is a disjoint union of  some $[x]$'s, where $[x]$ denotes the set of generators of $\langle x\rangle$.

\begin{prop}\label{identity}{\rm\cite[Proposition 4]{cam}}
Let $e$ be the identity of $G$.

{\rm(i)} If $G=\langle x\rangle$, then $\overline e$ is $G$ or $[e]\cup[x]$ according to $|x|$ is a prime power or not.

{\rm(ii)} If $G$ is a generalized quaternion  $2$-group, then $\overline e=[e]\cup[x]$, where $x$ is the unique involution in $G$.

{\rm(iii)} If $G$ is neither a cyclic group nor a generalized quaternion  $2$-group, then $\overline e=[e]$.
\end{prop}

\begin{prop}\label{closed}{\rm\cite[Proposition 5]{cam}}
Let $x$ be an element of $G$. Suppose $\overline x\neq\overline e$. Then one of the following holds.

{\rm(i)} $\overline x=[x]$.

{\rm(ii)} There exist distinct elements $x_1,x_2,\ldots,x_r$ in $G$ such that
$$
\overline x=[x_1]\cup[x_2]\cup\cdots\cup[x_r],\quad r\geq 2,
$$
where $\langle x_1\rangle\subseteq\langle x_2\rangle\subseteq\cdots\subseteq\langle x_r\rangle$,
$|x_i|=p^{s+i}$ for some prime $p$  and integer $s\geq 0$.
\end{prop}

The equivalence class $\overline e$ is said to be of {\em type} I.
An equivalence class  that does not contain $e$ is said to be of {\em type} II or  III according to  Proposition~\ref{closed} (i) or  (ii) holds. Furthermore, if $\overline x$ is of type III, with reference to Proposition~\ref{closed} (ii), the numbers $p, r, s$ are uniquely determined by $\overline x$.
We call $(p,r,s)$ its {\em parameters}.

 For each $x\in G$ and $\pi\in\aut(\mathcal P_G)$,  we have
$\overline x^\pi=\overline{x^\pi}.$ Hence,
  $\aut(\mathcal P_G)$ induces an action on $\mathcal U(G)$ as follow:
$$
\mathcal U(G)\times \aut(\mathcal P_G)\longrightarrow\mathcal U(G),\quad (\overline x,\pi)\longmapsto\overline{x^\pi}.
$$

Next we shall show that each orbit of $\aut(\mathcal P_G)$  on $\mathcal U(G)$ consists of some equivalence classes of the same type.

Note that $\overline e$ consists of vertices whose closed neighborhoods in $\mathcal P_G$ are $G$. Hence, one gets the following result.

\begin{lemma}\label{typeI}
  Each automorphism  of $\mathcal P_G$ fixes $\overline e$.
\end{lemma}

\begin{lemma}\label{un3}
 If $\overline x$ is an   equivalence class   of type III with parameters $(p,r,s)$, then  $|\overline x|=p^s(p^r-1)$.
\end{lemma}
\proof  With reference to Proposition~\ref{closed} (ii), we have
$$
|[x_i]|=\varphi(p^{s+i})=p^{s+i-1}(p-1),
$$
 which implies that
 $$
 |\overline x|=\sum_{i=1}^rp^{s+i-1}(p-1)=p^s(p^r-1),
 $$
  as desired.
$\qed$

\begin{lemma}\label{un4}
   Suppose $\overline x$ and $\overline y$ are two distinct  equivalence classes    of  type II or III. If $\langle x\rangle\subset\langle y\rangle$, then $|\overline x|\leq|\overline y|$, with equality  if and only if the follows hold.

  {\rm(i)} Both $\overline x$ and $\overline y$ are of type II.

  {\rm(ii)} $|y|=2|x|$ and $|x|$ is odd at least $3$.
\end{lemma}
\proof We divide the proof in three cases:

\medskip
{\em Case 1.} $\overline x$ is of type III with parameters $(p,r,s)$.

 Pick elements $x_1$ and $x_r$ in $\overline x$ of order $p^{s+1}$ and $p^{s+r}$, respectively. Then $\langle x_1\rangle\subseteq\langle x\rangle\subseteq\langle x_r\rangle$. Since $\langle x\rangle\subset\langle y\rangle$, we have $y\in N[x]=N[x_r]$. Note that any element $z$ satisfying $\langle x\rangle\subseteq\langle z\rangle\subseteq\langle x_r\rangle$ belongs to $\overline x$. Then $\langle x_r\rangle\subset\langle y\rangle$. Since $\frac{|y|}{p^{s+r}}$ is more than $1$, it has a prime divisor $p'$. Pick an element $z_0$ in $\langle y\rangle$ of order $p^{s+1}p'$.  Then $z_0\in N[x_1]=N[x_r]$. Hence, one of $p^{s+1}p'$ and $p^{s+r}$ is divided by the other. In view of $r\geq 2$, we get $p'=p$. It follows that $p^{s+r+1}$ divides $|y|$, and so $|[y]|\geq p^{s+r}(p-1)\geq p^{s+r}$.  Lemma~\ref{un3} implies that $|\overline x|<p^{s+r}$. Because $[y]\subseteq\overline y$, one has $|\overline x|<|\overline y|$.

\medskip
{\em Case 2.} $\overline y$ is of type III with parameters $(q,t,j)$.

Pick an element $y_1$   in $\overline y$ of order $q^{j+1}$. Since any element $z$ satisfying
$\langle y_1\rangle\subseteq\langle z\rangle\subseteq \langle y\rangle$ belongs to $\overline y$, one gets $\langle x\rangle\subset \langle y_1\rangle$.
Pick $y_0\in\langle y_1\rangle$ of order $q^j$. Then $\overline x\subseteq\langle y_0\rangle\setminus\{e\}$.
Hence $|\overline x|\leq q^j-1<q^j(q^t-1)$. According to Lemma~\ref{un3}, we have $|\overline x|<|\overline y|$.

\medskip
{\em Case 3.}  $\overline x$ and $\overline y$ are of type II.

 Then $|\overline x|=\varphi(|x|)$ and $|\overline y|=\varphi(|y|)$. Since $|x|$ divides $|y|$, it follows  that $|\overline x|$ divides $|\overline y|$, and so $|\overline x|\leq|\overline y|$. Note that $|x|\neq|y|$ and $x\neq e$. Then $|\overline y|=|\overline x|$ if and only if $|y|=2|x|$ and $|x|$ is odd at least $3$.

\medskip
Combining all these cases, we get the desired result.
$\qed$

\begin{lemma}\label{un1}
 Suppose $\overline x$ and $\overline y$ are two distinct  equivalence classes   of  type II or III. If $\langle x\rangle\subset\langle y\rangle$, then $\langle x^\pi\rangle\subset\langle y^\pi\rangle$ for every automorphism $\pi\in\aut(\mathcal P_G)$.
\end{lemma}
\proof Denote by $E_G$    the edge set  of $\mathcal P_G$.
Since $\{x,y\}\in E_G$, one gets $\{x^\pi,y^\pi\}\in E_G$. Because $\overline x^\pi\neq\overline y^\pi$, we have $\langle x^\pi\rangle\subset\langle y^\pi\rangle$ or $\langle y^\pi\rangle\subset\langle x^\pi\rangle$. Suppose for the contrary that $\langle y^\pi\rangle\subset\langle x^\pi\rangle$. By Lemma~\ref{un4}, we have $|\overline x|\leq|\overline y|$ and $|\overline y^\pi|\leq|\overline x^\pi|$. The fact that $\pi$ is a bijection implies that $|\overline x|=|\overline x^\pi|=|\overline y|=|\overline y^\pi|$. By Lemma~\ref{un4} again, the following hold:

a) For each $u\in\{x,y,x^\pi,y^\pi\}$,   $\overline u$ is of type II.

b) $|y|=|x^\pi|=2|x|=2|y^\pi|$ and $|y^\pi|$ is odd at least $3$.

Pick an element $z$ of order $2$ in $\langle y\rangle$. Then $\{z,y\}\in E_G$ and $\{z,x\}\not\in E_G$, which imply that $\{z^\pi,y^\pi\}\in E_G$ and $\{z^\pi,x^\pi\}\not\in E_G$, and hence $\langle y^\pi\rangle\subset\langle z^\pi\rangle$. Consequently,
$$
|\overline z^\pi|=|\overline{z^\pi}|\geq|[z^\pi]|=\varphi(|z^\pi|)\geq\varphi(|y^\pi|)\geq 2.
$$
Since $\overline z$ is of type II,  we get $|\overline z|=\varphi(2)=1$, a contradiction.   $\qed$

\begin{lemma}\label{un2}
  Suppose    $\overline x$ is of type II or III.     If   $|x|$ is a power of  a prime $p$, then $|x^\pi|$ is also a power of $p$ for any $\pi\in\aut(\mathcal P_G)$.
\end{lemma}
\proof Pick any prime divisor $q$ of $|x^\pi|$. It suffices to prove $q=p$. We only need to consider that $G$ is not a $p$-group. Let $z^\pi$ be an element of order $q$ in $\langle x^\pi\rangle$.  Proposition~\ref{identity} implies that $\overline z^\pi$ is of type II or III. It follows from Lemma~\ref{typeI} that $\overline z$ is of type II or III.

\medskip
{\em Claim 1. $|\overline z^\pi|=q^r-1$ for some positive integer $r$.}

If $\overline z^\pi$ is of type II, then $|\overline z^\pi|=\varphi(q)=q-1$. If $\overline z^\pi$ is of type III, then its parameters are $(q,r,0)$, which implies that $|\overline z^\pi|=q^r-1$ by Lemma~\ref{un3}.

\medskip
{\em Claim 2. $|\overline z|=p^j-1$ for some positive integer $j$.}

Since $\{z^\pi,x^\pi\}\in E_G$, one gets $\langle z\rangle\subseteq\langle x\rangle$ or $\langle x\rangle\subseteq\langle z\rangle$.
The fact that $z\neq e$   implies that $p$ divides $|z|$. Pick an element $y\in\langle z\rangle$ of order $p$. Note that $\overline y$ is  of type II or III. Similar to the proof of Claim 1, we get $|\overline y|=p^j-1$ for some positive integer $j$. It suffices to show that $\overline y=\overline{z}$. Suppose for the contrary that $\overline y\neq\overline{z}$. Then $\langle y\rangle\subset\langle z\rangle$. It follows from Lemma~\ref{un1} that $\langle y^{\pi}\rangle\subset\langle z^\pi\rangle$. Since   $|z^\pi|$ is a prime,  one has $y^{\pi}=e$. It follows that $\overline y^{\pi}$ is of type I,  contrary to Lemma~\ref{typeI}.

\medskip

Combining Claims 1 and 2, we get $q^r-1=p^j-1$, and so $q=p$, as desired.
$\qed$

\begin{prop}\label{un5}
  Let $\overline x\in\mathcal U(G)$ and  $\pi\in\aut(\mathcal P_G)$. Then $\overline x$ and $\overline x^\pi$ are of the same type. Moreover, if $\overline x$ is of type III, then $\overline x^\pi$ and $\overline x$ have the same parameters.
\end{prop}
\proof Suppose that $\overline x$ and $\overline x^\pi$ are of the distinct types. From Lemma~\ref{typeI}, we may assume that $\overline x$ is of type II and $\overline x^\pi$ is of type III with parameters $(p,r,s)$. Then $|\overline x^\pi|=p^s(p^r-1)$ by Lemma~\ref{un3}. Since  $|x^\pi|$ is a power of $p$, it follows from Lemma~\ref{un2} that $|x|=p^m$ for some positive integer $m$. Then $|\overline x|=|[x]|=\varphi(p^m)=p^{m-1}(p-1)$. Consequently, we get $p^s(p^r-1)=p^{m-1}(p-1)$, and so $r=1$, a contradiction. Therefore $\overline x$ and $\overline x^\pi$ are of the same type.

Suppose $\overline x$ and $\overline x^\pi$ are of type III with parameters $(p_1,r_1,s_1)$ and $(p_2,r_2,s_2)$, respectively. According to Lemmas~\ref{un3} and \ref{un2}, we get $p_1^{s_1}(p_1^{r_1}-1)=p_2^{s_2}(p_2^{r_2}-1)$ and $p_1=p_2$, and so $(p_1,r_1,s_1)=(p_2,r_2,s_2)$, as desired.
$\qed$

\section{Proof of main results}
In this section we present the proof of Theorems~\ref{mainthm1} and \ref{mainthm2}.
The following is an immediate result from (\ref{p(G)}), (\ref{s1}) and (\ref{s2}).

\begin{lemma} \label{lemma1}
Let $\pi$ be a permutation on the set $G$.

{\rm(i)} If $\pi\in P(G)$, then $\langle x\rangle^\pi=\langle x^\pi\rangle$ for each $x\in G$.

  {\rm(ii)} Then $\pi\in\prod_{i=1}^kS_{[C_i]}$ if and only if $[x]^\pi=[x]$ for each $x\in G$.

  {\rm(iii)} Then $\pi\in\prod_{i=1}^lS_{\overline{u_i}}$ if and only if $\overline x^\pi=\overline x$ for each $x\in G$.
\end{lemma}

\begin{lemma}\label{subgroups}
{\rm(i)} $P(G)$ and $\prod_{i=1}^kS_{[C_i]}$ are subgroups of $\aut(\overrightarrow{\mathcal P}_G)$.

{\rm(ii)} $P(G)$ and $\prod_{i=1}^lS_{\overline{u_i}}$ are subgroups of $\aut(\mathcal P_G)$.

\end{lemma}
\proof  (i) Pick $\sigma\in P(G)$ and $\xi\in\prod_{i=1}^kS_{[C_i]}$. In order to prove $\{\sigma,\xi\}\subseteq \aut(\overrightarrow{\mathcal P}_G)$, by (\ref{p(G)}) and (\ref{s1}), we only need to show that $(x,y)\in A_G$  implies $(x^\sigma,y^\sigma)\in A_G$  and $(x^\xi,y^\xi)\in A_G$, where   $A_G$ is the arc set of  $\overrightarrow{\mathcal P}_G$. Suppose $(x,y)\in A_G$. Then $\langle y\rangle\subseteq\langle x\rangle$.
It follows from Lemma~\ref{lemma1} that $\langle y^\sigma\rangle\subseteq\langle x^\sigma\rangle$ and $\langle y^\xi\rangle\subseteq\langle x^\xi\rangle$. Therefore $(x^\sigma,y^\sigma)\in A_G$
and $(x^\xi,y^\xi)\in A_G$.

(ii) Note that $\aut(\overrightarrow{\mathcal P}_G)\subseteq\aut(\mathcal P_G)$. By (i), we have $P(G)\subseteq\aut(\mathcal P_G)$.
Pick $\tau\in\prod_{i=1}^lS_{\overline{u_i}}$ and $\{x,y\}\in E_G$. By Lemma~\ref{lemma1}, we have $x^\tau\in\overline x$ and $y^\tau\in\overline y$.  If $\overline x=\overline y$, since $\overline x$ is a clique in $\mathcal P_G$, one has $\{x^\tau,y^\tau\}\in E_G$.
If $\overline x\neq\overline y$, then each vertex in $\overline x$ and each vertex in $\overline y$ are adjacent in $\mathcal P_G$, which implies that $\{x^\tau, y^\tau\}\in E_G$.
$\qed$

Write $\mathcal C'(G)=\{[x]\mid x\in G\}$. For each $[x]\in\mathcal C'(G)$ and $\pi\in\aut(\overrightarrow{\mathcal P}_G)$,  since
$[x]=\{x\}\cup\{y\mid \{(x,y),(y,x)\}\subseteq A_G\}$, we have
$$
[x]^\pi=\{x^\pi\}\cup\{y^\pi\mid\{(x^\pi,y^\pi),(y^\pi,x^\pi)\}\subseteq A_G\}=[x^\pi].
$$
Hence,  $\aut(\overrightarrow{\mathcal P}_G)$ induces an action on $\mathcal C'(G)$:
$$
\mathcal C'(G)\times \aut(\overrightarrow{\mathcal P}_G)\longrightarrow\mathcal C'(G),\quad ([x],\pi)\longmapsto[x^\pi].
$$

\begin{lemma}\label{normal}
{\rm(i)} $P(G)$ is a subgroup of the normalizer of $\prod_{i=1}^kS_{[C_i]}$ in $\aut(\overrightarrow{\mathcal P}_G)$.

{\rm(ii)} $P(G)$ is  a subgroup of the normalizer of $\prod_{i=1}^lS_{\overline{u_i}}$ in $\aut(\mathcal P_G)$.
\end{lemma}
\proof (i) Let $\sigma\in P(G)$ and $\xi\in\prod_{i=1}^kS_{[C_i]}$. For any $x\in G$, combining Lemmas~\ref{lemma1} and~\ref{subgroups}, we have
$$
[x]^{\sigma^{-1}\xi\sigma}=[x^{\sigma^{-1}}]^{\xi\sigma} =[x^{\sigma^{-1}}]^\sigma=[x].
$$
 It follows that $\sigma^{-1}\xi\sigma\in\prod_{i=1}^kS_{[C_i]}$, and so (i) holds.

(ii) The proof is similar to (i).
$\qed$

For each $\overline{u_i}\in\mathcal U(G)$, by Propositions~\ref{identity} and \ref{closed}, there exist  pairwise distinct $C_{i_1},\ldots,C_{i_t}\in\mathcal C(G)$ such that $\overline{u_i}=\bigcup_{j=1}^t[C_{i_j}]$. Hence, we get the following result.

\begin{lemma}\label{subgroup}
  $\prod_{i=1}^kS_{[C_i]}$ is a subgroup of $\prod_{i=1}^lS_{\overline{u_i}}$.
\end{lemma}

\begin{lemma}\label{intersection trivially}
 $|P(G)\cap(\prod_{i=1}^kS_{[C_i]})|=1$ and $|P(G)\cap(\prod_{i=1}^lS_{\overline{u_i}})|=1$.
 \end{lemma}
 \proof By Lemma~\ref{subgroup}, it is enough to prove $|P(G)\cap(\prod_{i=1}^lS_{\overline{u_i}})|=1$. Pick any $\pi\in P(G)\cap(\prod_{i=1}^lS_{\overline{u_i}})$ and  $x\in G$. Write $x=[C_i]_j$. Then $x^\pi=[C_i^\pi]_j$ by $\pi\in P(G)$. Since $\pi\in \prod_{i=1}^lS_{\overline{u_i}}$, one gets $x^\pi\in\overline x$.  Note that $\overline x$ is a clique in $\mathcal P_G$. Then $x^\pi\in N[x]$, and so $C_i\subseteq C_i^\pi$ or $C_i^\pi\subseteq C_i$.  Since $|C_i^\pi|=|C_i|$, we have $C_i^\pi=C_i$, which implies that $x^\pi=x$, as desired.
 $\qed$

For $x\in G$, we have $\langle x\rangle=\{x\}\cup\{y\mid (x,y)\in A_G\}$. Hence, for each $\pi\in\aut(\overrightarrow{\mathcal P}_G)$,
 $$
\langle x\rangle^\pi=\{x^\pi\}\cup\{y^\pi\mid (x^\pi,y^\pi)\in A_G\}=\langle x^\pi\rangle.
$$
Therefore,  $\aut(\overrightarrow{\mathcal P}_G)$ induces an action on $\mathcal C(G)$:
$$
\mathcal C(G)\times \aut(\overrightarrow{\mathcal P}_G)\longrightarrow\mathcal C(G),\quad (\langle x\rangle,\pi)\longmapsto\langle x^\pi\rangle.
$$
It is routine to verify that this group action preserves order, inclusion and noninclusion. Hence, the following result holds.

\begin{lemma}\label{same}
  For any $\pi\in\aut(\overrightarrow{\mathcal P}_G)$, there exists an element $\sigma\in P(G)$ such that $\langle x\rangle^\pi=\langle x\rangle^\sigma$ for every $x\in G$.
\end{lemma}

\noindent{\em Proof of Theorem~\ref{mainthm1}:}  It is apparent from Lemmas~\ref{normal} and \ref{intersection trivially} that $(\prod_{i=1}^k S_{[C_i]})\rtimes P(G)$ is a subgroup of $\aut(\overrightarrow{\mathcal P}_G)$.
Pick any $\pi\in\aut(\overrightarrow{\mathcal P}_G)$. By Lemma~\ref{same} there exists an element $\sigma\in P(G)$ such that, for any $x\in G,$
$$
\langle x^{\pi\sigma^{-1}}\rangle=\langle x\rangle^{\pi\sigma^{-1}}=\langle x\rangle,
$$
which implies that $x^{\pi\sigma^{-1}}\in[x]$. Then $\pi\sigma^{-1}\in\prod_{i=1}^k S_{[C_i]}$  and $\pi\in(\prod_{i=1}^k S_{[C_i]})(P(G))$. Hence, the desired result follows.
$\qed$

 \begin{prop}\label{un6}
 For any $\pi\in\aut(\mathcal P_G)$,  there exists an element $\tau\in\prod_{i=1}^lS_{\overline{u_i}}$ such that $\tau\pi\in\aut(\overrightarrow{\mathcal P}_G)$.
\end{prop}
\proof  Without loss of generality, assume that $\overline{u_1}$ is of type I, $\overline{u_i}$ is of type II for $2\leq i\leq d$, and $\overline{u_{j}}$ is of type III with parameters $(p_j,r_j,s_j)$ for $d+1\leq j\leq l$. According to Proposition~\ref{un5} each $\overline{u_j}^\pi$ is of type III with parameters $(p_j,r_j,s_j)$.

 For each $t\in\{1,\ldots,r_j\}$,  let $\{x_{t1}^{(j)},\ldots,x_{tm_{jt}}^{(j)}\}$ and $\{y_{t1}^{(j)},\ldots,y_{tm_{jt}}^{(j)}\}$ be the sets of elements of order $p_j^{s_j+t}$ in $\overline{u_j}$ and $\overline{u_j}^\pi$, respectively.
Then
$\tau_j: x_{tm}^{(j)}\longmapsto (y_{tm}^{(j)})^{\pi^{-1}}$
  is a permutation on $\overline{u_j}$, where $1\leq t\leq r_j$ and $1\leq m\leq m_{jt}$. Write $\tau=(\tau_1,\ldots,\tau_l)$, where $\tau_1$ is the inverse of the restriction of $\pi$ to $\overline{u_1}$, and $\tau_{i}$ is the identity of $S_{\overline{u_i}}$ for $2\leq i\leq d$. Hence $\tau\in\prod_{i=1}^l S_{\overline{u_i}}$.

We claim that,  for each $x\in G$, the equality  $|x^{\tau\pi}|=|x|$ holds. We divide our proof into three cases.

\medskip
 {\em Case 1.} $\overline x$ is of type I.

 Then $x^{\tau\pi}=(x^{\tau_1})^\pi=x$, and so $|x^{\tau\pi}|=|x|$.

\medskip
 {\em Case 2.} $\overline x$ is of type II.

 Then $x^{\tau\pi}=x^\pi$. According to Proposition~\ref{un5} we obtain that $\overline x^\pi$ is of type II, which implies that $[x^\pi]=\overline x^\pi=[x]^\pi$. Hence, one has
 $$
 \varphi(|x^\pi|)=|[x^\pi]|=|[x]|=\varphi(|x|).
$$
Suppose $|x^\pi|\neq|x|$. Without loss of generality, assume that $|x^\pi|<|x|$. Then $|x|=2|x^\pi|$ and $|x^\pi|$ is odd. Pick an element $z\in\langle x\rangle$ of order $2$. It is clear that $\overline z$ is of type II and $\overline z\neq\overline x$. From Lemma~\ref{un1} we get $\langle z^\pi\rangle\subset\langle x^\pi\rangle$. Since $\overline{z}^\pi$ is of type II, we infer that $|z^\pi|$ is odd at least $3$. Hence $|\overline z^\pi|\geq\varphi(3)=2$, contrary to $|\overline z|=1$. Therefore $|x^\pi|=|x|$.

\medskip
{\em Case 3.} $\overline x$ is of type III with parameters $(p_j,r_j,s_j)$.

Then  $x=x_{tm}^{(j)}$ for some indices $t$ and $m$. Since $x^{\tau\pi}=y_{tm}^{(j)}$, we have $|x^{\tau\pi}|=|x|$.

\medskip
Consequently, our claim is valid.

 Finally, we show that $\tau\pi\in\aut(\overrightarrow{\mathcal P}_G)$. Suppose $(u,v)\in A_G$. Then $\{u,v\}\in E_G$. It follows from Lemma~\ref{subgroups} that $\tau\pi\in\aut(\mathcal P_G)$. Hence $\langle u^{\tau\pi}\rangle\subseteq\langle v^{\tau\pi}\rangle$ or $\langle v^{\tau\pi}\rangle\subseteq\langle u^{\tau\pi}\rangle$.
 Since $\langle v\rangle\subseteq\langle u\rangle$, by the claim, $|v^{\tau\pi}|$ divides $|u^{\tau\pi}|$, which implies that  $\langle v^{\tau\pi}\rangle\subseteq\langle u^{\tau\pi}\rangle$.  So $(u^{\tau\pi},v^{\tau\pi})\in A_G$, as desired.
$\qed$

Combining  Theorem~\ref{mainthm1}, Lemmas~\ref{normal},  \ref{subgroup}, \ref{intersection trivially} and Proposition~\ref{un6}, we complete the proof of  Theorem~\ref{mainthm2}.

\section{Examples}
In this section we shall compute $\aut(\overrightarrow{\mathcal P}_G)$ and $\aut(\mathcal P_G)$ if $G$ is  cyclic,  elementary abelian,   dihedral or  generalized quaternion. We begin with   cyclic groups.

\begin{example}\label{zn}
  Let $n$ be a positive integer. Then

 {\rm(i)} $\aut(\overrightarrow{\mathcal P}_{Z_n})\cong\prod_{d\in D(n)}S_{\varphi(d)}.$

{\rm(ii)} $\aut(\mathcal P_{Z_n})\cong\left\{
  \begin{array}{ll}
  S_n,& \textup{if $n$ is a prime power},\\
S_{\varphi(n)+1}\times\prod_{d\in D(n)\setminus\{1,n\}}S_{\varphi(d)},&\textup{otherwise}.
  \end{array}\right.$
\end{example}
\proof For any $d\in D(n)$, denote by $A_d$ the unique cyclic subgroup of order $d$ in $Z_n$. Note that $P(Z_n)=\mathbf 1_{\{A_d\mid d\in D(n)\}}$ and $S_{[A_d]}\cong S_{\varphi(d)}$, where $\mathbf{1}_{\Omega}$ denotes the identity map on the set $\Omega$.
It follows from Theorem~\ref{mainthm1}  that (i) holds.
If $n$ is a prime power, then $\aut(\mathcal P_G)\cong S_n$ by \cite[Theorem 2.12]{cha}.
If $n$ is not a prime power,
by \cite[Proposition 3.6]{fmw},
$$
\mathcal U(Z_n)=\{[A_d]\mid d\in D(n)\setminus\{1,n\}\}\cup\{[A_1]\cup[A_n]\}.
$$
Hence (ii) holds by Theorem~\ref{mainthm2}.
$\qed$

Example~\ref{zn} shows that the conjecture   proposed by Doostabadi, Erfanian and   Jafarzadeh  holds if $n$ is not a prime power.

Combining Theorems~\ref{mainthm1} and ~\ref{mainthm2}, we get the following result.

\begin{prop}\label{sn}
  $\aut(\overrightarrow{\mathcal P}_G)=\aut(\mathcal P_G)$ if and only if $\overline x=[x]$ for each $x\in G$.
\end{prop}

%Let $N$ be a group and $H$ be a permutation group on $\Omega$. Define
%$$
%\kappa_h:\prod_{\omega\in\Omega}N_\omega\longrightarrow \prod_{\omega\in\Omega}N_\omega,\quad
%(n_{\omega_1},\ldots,n_{\omega_{|\Omega|}})\longmapsto
%(n_{\omega_1^{h^{-1}}},\ldots,n_{\omega_{|\Omega|}^{h^{-1}}}),
%$$
%where $N_w$ is the copy of $N$ and $n_{w_i}\in N$. Then the following is a homomorphism.
%$$
%\kappa: H\longrightarrow\prod_{\omega\in\Omega}N_\omega,\quad h\longmapsto\kappa_h.
%$$
%The wreath product $N\wr H$ is the semiproduct $\prod_{w\in\Omega}N_w\rtimes_\kappa H$.

Let $H$ be a group and $K$ be a permutation group on a set $Y$. The wreath product $H\wr K$ is the semidirect product $N\rtimes K$, where
$N$ is the direct product of $|Y |$ copies of $H$ (indexed by $Y$), and $K$ acts on $N$ by permuting the factors in the same way as it permutes elements of $Y$. If $H$ is a permutation group on a set $X$, then $H\wr K$ has a nature action on $X\times Y$:
$$
(X\times Y)\times (H\wr K)\longrightarrow X\times Y,\quad ((x,y_i),(h_{y_1},\ldots,h_{y_{|Y|}};k))\longmapsto(x^{h_{y_i}},y_i^k),
$$
where $\{y_1,\ldots,y_{|Y|}\}=Y$.

For a prime $p$ and a positive integer $n$, let $Z^n_p$ denote the elementary abelian $p$-group, i.e., the direct product of $n$ copies of $Z_p$.

\begin{example}
  Let $n\geq 2$.  Then
  $$
  \aut(\mathcal P_{Z_p^n})=\aut(\overrightarrow{\mathcal P}_{Z_p^n})\cong S_{p-1}\wr S_m,
  $$
  where $m=\frac{p^n-1}{p-1}$.
\end{example}
\proof Write $\mathcal C(Z^n_p)=\{\langle e\rangle, A_1,\ldots,A_m\}$. Then  each $A_i$ is isomorphic to $Z_p$ and $|A_i\cap A_j|=1$ for $i\neq j$. Hence, one has $P(Z^n_p)=\mathbf 1_{\{\langle e\rangle\}}S_{\{A_i\mid 1\leq i\leq m\}}$. Combining Theorems~\ref{mainthm1},   \ref{mainthm2} and Proposition~\ref{sn}, we get the desired result.
$\qed$

\begin{figure}[hptb]
\begin{center}
\includegraphics{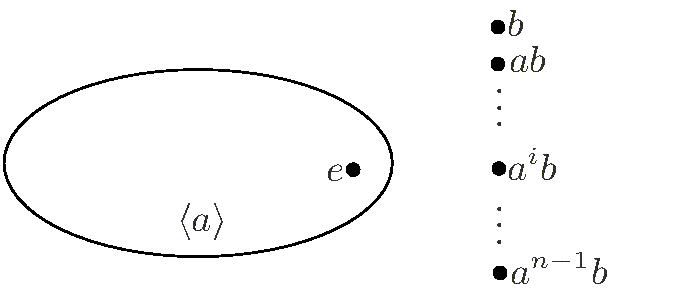}
\caption{The partition of $D_{2n}$\label{D2n} }
\end{center}
\end{figure}

\begin{example}
  For $n\geq 3$, let $D_{2n}$ denote the dihedral group of order $2n$. Then

 {\rm(i)}  $\aut(\overrightarrow{\mathcal P}_{D_{2n}})\cong \prod_{d\in D(n)}S_{\varphi(d)}\times S_n.$

 {\rm(ii)}  $\aut(\mathcal P_{D_{2n}})\cong \left\{
  \begin{array}{ll}
   S_{n-1}\times S_n,&\textup{if $n$ is a prime power},\\
   \prod_{d\in D(n)}S_{\varphi(d)}\times S_n,&\textup{otherwise}.
  \end{array}\right.$
\end{example}
\proof Pick $a, b\in D_{2n}$ with $|a|=n$ and $|b|=2$. Then
$$
\begin{array}{l}
D_{2n}=\{e,a,\ldots,a^{n-1}\}\cup\{b,ab,\ldots,a^{n-1}b\},\\
\mathcal C(D_{2n})=\mathcal C(\langle a\rangle)\cup\{\langle a^ib\rangle\mid 0\leq i\leq n-1\},
\end{array}
$$
as shown in Figure~\ref{D2n}.
Note that $|a^ib|=2$, $|\langle a^ib\rangle \cap\langle a\rangle|=1$ and $|\langle a^ib\rangle\cap\langle a^jb\rangle|=1$ for $i\neq j$.  Hence we have
\begin{equation}\label{d}
P(D_{2n})=\mathbf{1}_{\mathcal C(\langle a\rangle)} S_{\{\langle a^ib\rangle\mid 0\leq i\leq n-1\}}.
\end{equation}
By Theorem~\ref{mainthm1}, one has
$$
\aut(\overrightarrow{\mathcal P}_{D_{2n}})=(\prod_{A\in\mathcal C(\langle a\rangle)}S_{[A]}\rtimes \mathbf{1}_{\langle a\rangle})\times (\prod_{i=0}^{n-1}S_{[\langle a^ib\rangle]}\rtimes S_{\{ a^ib\mid 0\leq i\leq n-1\}}),
$$
 which implies (i).

Suppose $n$ is a prime power. Then
$$\mathcal U(D_{2n})=\{\{e\},\{a^1,\ldots,a^{n-1}\}\}\cup\{\{a^ib\}\mid 0\leq i\leq n-1\}.$$
Theorem~\ref{mainthm2} and  (\ref{d}) imply that $\aut(\mathcal P_{D_{2n}})\cong S_{n-1}\times S_n$.

Suppose $n$ is not a prime power. Then $\overline{a^i}=[a^i]$. By Proposition~\ref{sn}, our desired result follows.
$\qed$

\begin{figure}[hptb]
\begin{center}
\includegraphics{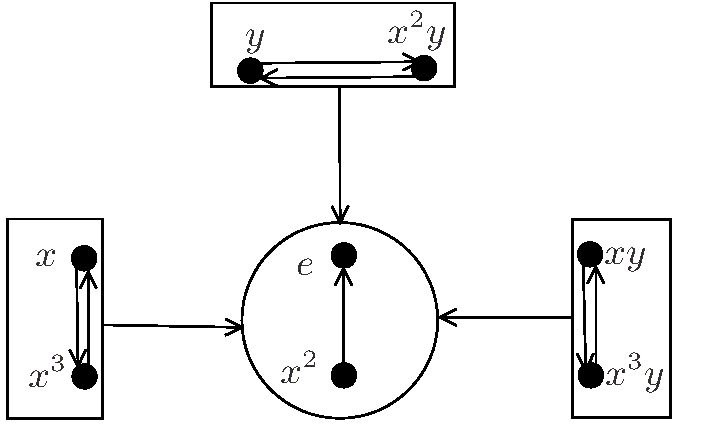}\caption{The power digraph of $Q_8$\label{Q8}}
\end{center}
\end{figure}

Let $Q_{4n}$ denote the generalized quaternion group of order $4n$, i.e.,
\begin{equation}\label{q4n}
Q_{4n}=\langle x,y\mid x^{2n}=e,x^n=y^2, y^{-1}xy=x^{-1}\rangle.
\end{equation}
The power digraph $\overrightarrow{\mathcal P}_{Q_8}$ is shown in Figure~\ref{Q8}. Observe that
  $$
  \aut(\overrightarrow{\mathcal P}_{Q_8})\cong  S_2\wr S_3,\quad \aut(\mathcal P_{Q_8})\cong S_2\times (S_2\wr S_3).
  $$

\begin{figure}[hptb]
\begin{center}
\includegraphics{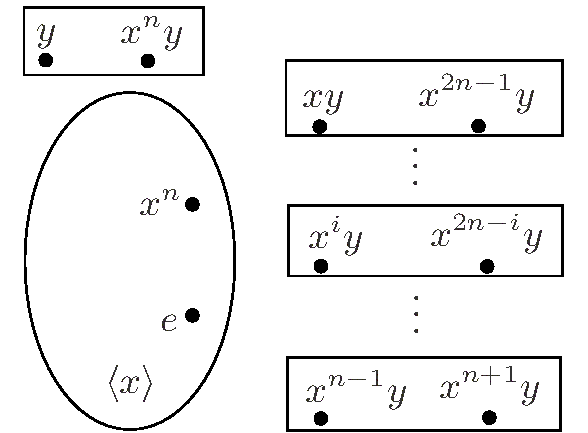}
\caption{The partition of $Q_{4n}$\label{Q4n}}
\end{center}
\end{figure}

\begin{example}
  Let $n\geq 3$. Then

{\rm(i)}   $\aut(\overrightarrow{\mathcal P}_{Q_{4n}})\cong\prod_{d\in D(2n)}S_{\varphi(d)}\times(S_2\wr S_n).$

{\rm(ii)} $\aut(\mathcal P_{Q_{4n}})\cong
  \begin{cases}
   S_2\times S_{2n-2}\times (S_2\wr S_n),&\textup{if $n$ is a power of 2},\\
   \prod_{d\in D(2n)}S_{\varphi(d)}\times(S_2\wr S_n),&\textup{otherwise}.
  \end{cases}
 $
\end{example}
\proof With reference to $(\ref{q4n})$, $y^{-1}=x^ny$ and $(x^iy)^{-1}=x^{2n-i}y$ for $i\in\{1,\ldots,n-1\}$. So we have  $$
\begin{array}{l}
Q_{4n}=\{e,x,\ldots,x^{2n-1}\}\cup\{y,x^ny\}\cup\bigcup_{i=1}^{n-1}\{x^iy,x^{2n-i}y\},\\
\mathcal C(Q_{4n})=\mathcal C(\langle x\rangle)\cup\{\langle x^jy\rangle\mid 0\leq j\leq n-1\},
\end{array}
$$
as shown in Figure~\ref{Q4n}.
Then
\begin{equation}\label{q}
P(Q_{4n})=\mathbf 1_{\mathcal C(\langle x\rangle)} S_{\{\langle x^jy\rangle\mid 0\leq j\leq n-1\}}.
\end{equation}
Thus (i) holds from Theorem~\ref{mainthm1}.

Suppose $n$ is a power of $2$. Then
$$
\mathcal U(Q_{4n})=\{\{e,x^n\},\langle x\rangle\setminus\{e, x^n\},\{y,x^ny\}\}\cup\{\{x^iy,x^{2n-i}y\}\mid 1\leq i\leq n-1\}.
$$
Theorem~\ref{mainthm2} and (\ref{q}) imply (ii) holds.

Suppose $n$ is not a power of $2$. Then $\overline{x^i}=[x^i]$ for each $i\in\{0,1,\ldots,2n-1\}$. From Proposition~\ref{sn}  we get the desired result.
$\qed$

\section*{Acknowledgement}

This research is supported
by National Natural Science Foundation of China (11271047,  11371204).

\end{CJK*}

\end{document}